\documentclass{article}
\usepackage{amsmath}
\usepackage{amsthm}
\usepackage{amsfonts}
\usepackage{amssymb}
\usepackage{enumerate}
\usepackage{comment}
\usepackage{hyperref}

\newtheorem{theorem}{Theorem}
\newtheorem{lemma}{Lemma}

\theoremstyle{remark}
\newtheorem{remark}{Remark}

\title{Asymptotic Stabilizability of Underactuated Hamiltonian Systems With Two Degrees of Freedom\thanks{S.D. Grillo and L.M. Salomone thank CONICET for its financial support.}}

\author{Sergio D. Grillo\thanks{Instituto Balseiro, UNCuyo-CNEA, San Carlos de Bariloche, R\'{\i}o Negro, Rep\'{u}blica Argentina  (\href{mailto:sergiog@cab.cnea.gov.ar}{sergiog@cab.cnea.gov.ar}).}
\and Leandro M. Salomone\thanks{CMaLP, Fac. de Ciencias Exactas, UNLP, La Plata, Buenos Aires, Rep\'{u}blica Argentina
  (\href{mailto:salomone@mate.unlp.edu.ar}{salomone@mate.unlp.edu.ar}, \href{mailto:marcezuccalli@gmail.com}{marcezuccalli@gmail.com}).}
\and Marcela Zuccalli\footnotemark[3]}

\begin{document}

\maketitle

\begin{abstract}
For an underactuated (simple) Hamiltonian system with two degrees
of freedom and one degree of underactuation, a rather general condition that ensures its stabilizability,
by means of the existence of a (simple) Lyapunov function, was found
in a recent paper by D.E. Chang within the context of the energy shaping
method. Also, in the same paper, some additional assumptions were
presented in order to ensure also asymptotic stabilizability. In this
paper we extend these results by showing that above mentioned condition
is not only sufficient, but also a necessary one. And, more importantly,
we show that no additional assumption is needed to ensure asymptotic
stabilizability. 
\end{abstract}


\section{Introduction}
Consider an underactuated Hamiltonian system with
two degrees of freedom and exactly one actuator (i.e. with one degree of underactuation). Such a system can
be described by a pair $\left(H,Y\right)$, where $H$ is a Hamiltonian
function on a $4$-dimensional phase space and $Y$ is a (vertical)
vector field defining the direction of the actuator. Fix a critical
point $\alpha_{0}$ for the related Hamiltonian vector field $X_{H}$
and assume from now on that:
\begin{enumerate}
\item $H$ is \textit{simple}, i.e. in any canonical coordinate chart $\left(x,y,p_{x},p_{y}\right)$
\begin{equation*}
H\left(x,y,p_{x},p_{y}\right) =\frac{1}{2}\,\left(p_{x},p_{y}\right)\begin{bmatrix}a(x,y) & b(x,y)\\
b(x,y) & c(x,y)
\end{bmatrix}\,\left(\begin{array}{l}
p_{x}\\
p_{y}
\end{array}\right) +h\left(x,y\right)
\end{equation*}
for some function $h$ and with
\begin{equation*}
\mathbb{H}(x,y):=\begin{bmatrix}a(x,y) & b(x,y)\\
b(x,y) & c(x,y)
\end{bmatrix}\label{HV-1}
\end{equation*}
a positive-definite matrix;
\item there exist canonical coordinates, that we shall call \textbf{adapted
coordinates}, in which:
\begin{enumerate}
\item $Y$ is given by the constant vector $\left(0,0,0,1\right)$;
\item the critical point $\alpha_{0}$ is represented by the vector $\mathbf{\bar{0}}:=\left(0,0,0,0\right)$. 
\end{enumerate}
\end{enumerate}
A system satisfying above conditions will be called \textbf{underactuated
simple Hamiltonian system with two degrees of freedom}. In Ref. \cite{chang2},
D.E. Chang found, among other things, a sufficient condition that ensures the stabilizability\footnote{By \textbf{\emph{stabilizable}} at a point $\alpha_{0}$ we mean that
there exists a state feedback controller such that the related closed-loop
system is stable at $\alpha_{0}$.} of such systems 
at the given critical point $\mathbf{\bar{0}}$. The mentioned condition
can be written (in adapted coordinates) as
\begin{equation}
\left(b\,\frac{\partial^{2}h}{\partial x^{2}}+c\,\frac{\partial^{2}h}{\partial x\partial y}\right)\left(0,0\right)\neq0\ \ \ \text{or\ \ \ }\frac{\partial^{2}h}{\partial x^{2}}\left(0,0\right)>0.\label{mc}
\end{equation}
In coordinate-free terms, according to \cite{chang2}, above inequalities
mean that: either the linearization of the system at the given critical
point is controllable, or it is uncontrollable with uncontrollable
modes given by a purely imaginary pair.

\bigskip{}

Chang's work was done within the framework of the \emph{energy shaping}
method (see for instance \cite{auck,ks,cl,es1,es2,CBLMW02,kri,idapbc,shaft2,woosley}),
or more precisely, within his version of the method, developed
in \cite{chang1,chang3,chang4}. Let us briefly review such a method
in the present context. Its main idea is to construct, for a given
pair $(H,Y)$ and a given critical point of $X_{H}$, a state feedback
controller $u$ and a simple Lyapunov function $\hat{H}$ for the
resulting closed-loop system. Note that $\hat{H}$ (to be simple)
must have the form 
\begin{align}
\hat{H}\left(x,y,p_{x},p_{y}\right)= \frac{1}{2}\,\left(p_{x},p_{y}\right)\begin{bmatrix}f(x,y) & g(x,y)\\
g(x,y) & l(x,y)
\end{bmatrix}\,\left(\begin{array}{l}
p_{x}\nonumber\\
p_{y}
\end{array}\right)+\hat{h}\left(x,y\right)\label{V}
\end{align}
with 
\begin{equation}
\hat{\mathbb{H}}(x,y):=\begin{bmatrix}f(x,y) & g(x,y)\\
g(x,y) & \,l(x,y)
\end{bmatrix}\label{HV1-1}
\end{equation}
positive-definite. To find the controller $u$, a set of partial differential
equations (PDEs), known as \emph{matching conditions}, must be solved.
Such PDEs have the pair $(H,Y)$ as datum and the aforementioned Lyapunov
function $\hat{H}$ as their unknown. According to Ref. \cite{lyap-shap},
for pairs $\left(H,Y\right)$ as described above, and using adapted
coordinates, the mentioned PDEs are
\begin{equation}
\sum_{i,j,k,l=1}^{2}\left(\partial_{k}\hat{\mathbb{H}}^{ij}\,\mathbb{H}^{kl}-\partial_{k}\mathbb{H}^{ij}\,\hat{\mathbb{H}}^{kl}\right)\,p_{i}p_{j}p_{l}=0,\label{kec}
\end{equation}
the \emph{kinetic matching condition}, and 
\begin{equation}
\sum_{k,l=1}^{2}\left(\partial_{k}\hat{h}\,\mathbb{H}^{kl}-\partial_{k}h\,\hat{\mathbb{H}}^{kl}\right)\,p_{l}=0,\label{pec}
\end{equation}
the \emph{potential matching condition}, and must be satisfied for
all $(x,y,p_{x},p_{y})$ such that 
\begin{equation}
p_{x}\,g(x,y)+p_{y}\,l(x,y)=0.\label{lsur}
\end{equation}
Here, $\partial_{1}$ (resp. $\partial_{2}$) denotes the partial
derivative w.r.t. $x$ (resp. $y$), $p_{1}=p_{x}$ and $p_{2}=p_{y}$.
Note that for $\hat{H}$ to be a Lyapunov function related to the
point $\mathbf{\bar{0}}$, we need that $\hat{h}$ is a positive-definite
function w.r.t. $\left(0,0\right)$, i.e.
\begin{equation}
\hat{h}\left(0,0\right)=0\ \ \textrm{and}\ \ \hat{h}\left(x,y\right)>0\ \ \text{for all\ \ }\left(x,y\right)\neq\left(0,0\right).\label{pc}
\end{equation}

Once a solution $\hat{H}$ of the matching conditions is given, the
method provide a concrete procedure to construct a state feedback
controller. In the case under
consideration, such a controller has the form (see Ref. \cite{lyap-shap})
\begin{equation}
u(x,y,p_{x},p_{y})=\left(0,0,0,\lambda\left(x,y,p_{x},p_{y}\right)\right),\label{Y}
\end{equation}
with 
\begin{equation}
\lambda(x,y,p_{x},p_{y}):=-\frac{\left(\mu+\{\hat{H},H\}\right)(x,y,p_{x},p_{y})}{p_{x}\,g(x,y)+p_{y}\,l(x,y)}.\label{sol}
\end{equation}
Here $\{\cdot,\cdot\}$ denotes the canonical Poisson bracket and
$\mu$ is any non-negative function such that 
\begin{equation}
\frac{\mu(x,y,p_{x},p_{y})}{p_{x}\,g(x,y)+p_{y}\,l(x,y)}\label{muci}
\end{equation}
is smooth. Thus, if a solution $\hat{H}$ of the matching conditions \eqref{kec} \eqref{pec} and \eqref{lsur}
is found, and satisfies above mentioned positivity requirements {[}see
\eqref{pc}{]}, the system in question can be stabilized at the
point $\mathbf{\bar{0}}$ by means of the controller \eqref{Y}. 
And, as we said above, such stability is ensured by the existence of a Lyapunov function 
for the related closed-loop system: the solution $\hat{H}$.

\bigskip{}

It was shown in \cite{lyap-shap} that the Chang's version of the energy shaping method \cite{chang1,chang3,chang4} is \emph{maximal} among the so-called ``Lyapunov based methods.'' More precisely, if an underactuated simple Hamiltonian system (with any number of degrees of freedom and any degree of underactuation) is stabilized by a method that gives rise to a closed-loop system with a simple Lyapunov function, then such a function must be a solution of the matching conditions and the related controller is exactly the one given by the energy shaping method.

\bigskip{}

Coming back to \eqref{mc}, what Chang showed in \cite{chang2} was actually that \eqref{mc} is a necessary and sufficient condition for finding a simple solution of the matching conditions (with the above mentioned positive requeriments), and consequently a sufficient condition for stabilizability. Moreover, in the same paper, two additional assumptions to ensure not only stabilizability, but also asymptotic stabilizability, were presented.

\bigskip{}

In the present paper, we further study the stabilizability condition \eqref{mc}
and show:

\begin{description}
\item [$\mathbf{a}$.] A (slightly) different stabilizability characterization: \eqref{mc} is a sufficient and a necessary condition 
to stabilize an underactuated simple Hamiltonian system with
two degrees of freedom by any method (not only the energy shaping) that guarantees such stability by exhibiting
(or at least by ensuring the existence of) a simple Lyapunov function.
\item [$\mathbf{b}$.] The main result of the paper: \eqref{mc} not only ensures stabilizability,
but also asymptotic stabilizability. That is to say, no additional
condition is needed, other than \eqref{mc}, in order to prove the
asymptotic stabilizability for an underactuated simple Hamiltonian
system with two degrees of freedom. 
\end{description}

\bigskip{}

The paper is organized as follows. In \S\ref{bidimensional} we write down a more convenient expression of
the matching conditions \eqref{kec}, \eqref{pec} and \eqref{lsur}.
Then, studying the existence of their solutions, we give an alternative
derivation of the condition \eqref{mc} and show the point \textbf{$\mathbf{a}$} above. We could show that point simply by combining 
the results of Chang in \cite{chang2} and the above mentioned \emph{maximal character} of the energy
shaping method (showed in Ref. \cite{lyap-shap}). Nevertheless, we decided to make a detailed proof 
because of the involved reasoning and calculations, which are necessary to prove the second result of the paper. The latter is done in 
\S\ref{asymptotic}, where, by combining the LaSalle's invariance principle, a Dirac-like algorithm and the Morse Lemma, we show that condition \eqref{mc}
also implies asymptotic stabilizability (i.e. we prove the claim of
point $\mathbf{b}$). Finally, we illustrate our results with an
example.

\section{Stabilizability of systems with two degrees of freedom}\label{bidimensional} 
In this section we prove that, given an underactuated simple Hamiltonian system with two degrees of freedom, condition \eqref{mc} ensures its stabilizability. Reciprocally, if such a system is stabilizable and that stability can be ensured by a simple Lyapunov function, we prove that condition \eqref{mc} must hold. All that will be done by
using the Chang version of the energy shaping method \cite{chang1,chang3,chang4}. To begin with, we shall write down \eqref{kec}, \eqref{pec} and \eqref{lsur} in a way which is more appropriate for our purposes.

\bigskip{}

We want to emphasize that the results of the present section represent a slightly modification of those contained in the work of Chang  \cite{chang2} and, as we said in the Introduction, they can be proved just by using the maximal character of the energy shaping method \cite{lyap-shap}. However, we decide to give an alternative proof here because some of the intermediate steps are crucial for showing the main result of this paper, developed in \S\ref{asymptotic}.

\subsection{Re-writing the matching conditions}

Consider a pair $\left(H,Y\right)$ and a set of adapted coordinates $\left(x,y,p_{x},p_{y}\right)$
as those described in the Introduction. 
\begin{remark}
\label{critt}Since $H$ is simple, it is easy to show that $\overline{\mathbf{0}}$
is critical for the Hamiltonian vector field 
\begin{equation}
X_{H}=\left(\frac{\partial H}{\partial p_{x}},\frac{\partial H}{\partial p_{y}},-\frac{\partial H}{\partial x},-\frac{\partial H}{\partial y}\right)\label{XH}
\end{equation}
if and only if $\left(0,0\right)$ is critical for the function $h$.
\end{remark}
Consider also the matching conditions \eqref{kec}, \eqref{pec} and \eqref{lsur} for the unknowns $\hat{\mathbb{H}}$ and $\hat{h}$.
Note first that, since $\mathbb{H}$ is positive-definite (see 
\eqref{HV-1} in the Introduction), then 
\begin{equation}
a,c>0\;\;\;\textrm{and}\;\;\;\Delta:=a\,c-b^{2}>0.\label{pd}
\end{equation}
 Analogously, regarding $\hat{\mathbb{H}}$, we have that 
\begin{equation}
f>0\quad\text{and}\quad f\,l-g^{2}>0,\label{sc}
\end{equation}
which implies that $l>0$. To further simplify the notation, define
\begin{equation}
\delta:=\frac{f\,l-g^{2}}{\,l},\quad\gamma:=\frac{g}{\,l},\label{dl}
\end{equation}
and 
\begin{equation}
B:=a_{x}-2\,\gamma\,b_{x}+\gamma^{2}\,c_{x}.\label{AB}
\end{equation}
(From now on, the subindices $x$ and $y$ denote partial differentiation).
In terms of these new variables, it can be shown that the kinetic
matching condition \eqref{kec}, combined with \eqref{lsur}, and
the positivity conditions \eqref{sc} are equivalent to 
\begin{equation}
(a-b\,\gamma)\,\delta_{x}+(b-c\,\gamma)\,\delta_{y}=B\,\delta,\ \ \ \delta>0\label{regkinetic2}
\end{equation}
and $l>0$, while the potential matching condition \eqref{pec}, combined
with \eqref{lsur}, adopt the form 
\begin{equation}
(a-b\,\gamma)\,\hat{h}_{x}+(b-c\,\gamma)\,\hat{h}_{y}=h_{x}\,\delta.\label{regpotential2}
\end{equation}

Summing up, the matching and positivity conditions can be described
by \eqref{regkinetic2}, \eqref{regpotential2} and \eqref{pc}
for the unknowns $\left(\delta,\gamma,\hat{h}\right)$, plus the condition
$l>0$. To go back to the original variables, we just must use the
formulae {[}see \eqref{dl}{]}
\begin{equation}
f=\delta+l\,\gamma^{2}\;\;\;\textrm{and}\;\;\;g=l\,\gamma.\label{invv}
\end{equation}

\subsection{A sufficient condition for stabilizability}

\label{exandun}

As we said in the Introduction about the energy shaping method, if we find a solution to \eqref{kec}, \eqref{pec}, \eqref{lsur} and \eqref{pc}, then we can construct
a vector field $u$ {[}see \eqref{Y} and \eqref{sol}{]} and
a simple Lyapunov function that ensures the stability of the related
closed-loop system at $\mathbf{\bar{0}}$ (at least locally around
$\mathbf{\bar{0}}$). So, according to the last subsection, the
stabilizability of $\left(H,Y\right)$ around $\mathbf{\bar{0}}$
can be analyzed by studying the existence of solutions $\left(\delta,\gamma,\hat{h}\right)$
of \eqref{regkinetic2}, \eqref{regpotential2} and \eqref{pc}. To
do that, let us consider the next two lemmas. In what follows,
we shall call $U\subset\mathbb{R}^{4}$ the neighborhood of $\overline{\mathbf{0}}$
where the adapted coordinates take their values. Also, for simplicity,
we shall write $\mathbf{0}:=\left(0,0\right)$. 

\begin{lemma}\label{lem:gamma-1} Given a function $\gamma$ satisfying
\begin{equation}
\gamma\left(\mathbf{0}\right)\neq\frac{b\left(\mathbf{0}\right)}{c\left(\mathbf{0}\right)}\label{gbc}
\end{equation}
and 
\begin{equation}
\left[(a-b\,\gamma)\,h_{xx}+(b-c\,\gamma)\,h_{xy}\right]\left(\mathbf{0}\right)>0,\label{condpos}
\end{equation}
there exist functions $\delta$ and $\hat{h}$ such that $\left(\delta,\gamma,\hat{h}\right)$
is a solution of \eqref{regkinetic2}, \eqref{regpotential2} and
\eqref{pc}. \end{lemma} 
\begin{proof}
Let us begin with \eqref{regkinetic2}. This is a first-order PDE,
so we can use the \emph{Method of Characteristics} to find a solution
around $\mathbf{0}$. But, in order for this to make sense, we need
a suitable boundary condition on a non-characteristic submanifold
$\Gamma$. Let $V$ be the projection of $U$
onto the first two coordinates. Observe that the characteristic vector
field is $A=\left(a-b\gamma,b-c\gamma\right)$. Then, we may take
the submanifold $\Gamma\subset\mathbb{R}^{2}$ to be the $x$-axis,
i.e. to take 
\[
\Gamma=\left\{ \left(x,0\right):x\in\mathbb{R}\right\} \cap V,
\]
so long as we ensure that the second component of $A$ is nonzero
around $\mathbf{0}$. But this amounts to choose $\gamma$ such that
\eqref{gbc} holds. Since we need $\delta>0$, we can impose
the boundary condition $\left.\delta\right\vert _{\Gamma}=\,s$, where
$\,s:\mathbb{R}\rightarrow\mathbb{R}$ is a function such that $\,s(0)>0$.
In this case, the Theorem of Characteristics states that there is
a unique solution $\delta$ of \eqref{regkinetic2} such that $\delta(x,0)=\,s(x)$,
which implies, by continuity, that $\delta(x,y)>0$ around $\mathbf{0}$.
We can shrink $V$ (together with $U$), if necessary, in order to ensure that $\delta>0$
along all of $V$. From now on, we shall use this shrinking process
implicitly (finitely many times).

Let us continue with \eqref{regpotential2} and \eqref{pc}. The former
is also a first-order PDE, and with the same characteristic vector
field $A$. Assuming \eqref{gbc} again, the $x$-axis is a non-characteristic
submanifold and we can impose $\left.\hat{h}\right\vert _{\Gamma}=\,r$,
where $\,r:\mathbb{R}\rightarrow\mathbb{R}$ is a smooth function
such that $\,r(0)=0$. This implies that $\hat{h}\left(\mathbf{0}\right)=0$,
which is the first part of \eqref{pc}. The second part says that
$\mathbf{0}$ is an isolated minimum for $\hat{h}$, or equivalently,
$\mathbf{0}$ is critical for $\hat{h}$ and the Hessian of $\hat{h}$
is positive-definite at $\mathbf{0}$. Let us analyze these conditions.
Since $\mathbf{0}$ is critical for $h$ (see Remark \ref{critt}),
i.e. $\left(h_{x}(\mathbf{0}),h_{y}(\mathbf{0})\right)=\mathbf{0}$,
it follows from \eqref{regpotential2} that $\hat{h}$ must satisfy
\[
\left[(a-b\,\gamma)\,\hat{h}_{x}+(b-c\,\gamma)\,\hat{h}_{y}\right]\left(\mathbf{0}\right)=0.
\]
Thus, since $b(\mathbf{0})-c(\mathbf{0})\,\gamma(\mathbf{0})\neq0$
{[}see \eqref{gbc}{]}, in order to have that $\left(\hat{h}_{x}(\mathbf{0}),\hat{h}_{y}(\mathbf{0})\right)=\mathbf{0}$,
it suffices to ask that $\hat{h}_{x}(\mathbf{0})=\,r^{\prime}(0)=0$.
So far, we have that $r$ and $s$ must satisfy 
\begin{equation}
\,s\left(0\right)>0,\ \ r\left(0\right)=\,r^{\prime}(0)=0.\label{pr1}
\end{equation}
On the other hand, the Hessian of $\hat{h}$ is positive-definite
at $\mathbf{0}$ if and only if 
\begin{equation}
\hat{h}_{xx}(\mathbf{0})>0\qquad\text{and}\qquad\left(\hat{h}_{xx}\,\hat{h}_{yy}-\hat{h}_{xy}^{2}\right)(\mathbf{0})>0.\label{hes}
\end{equation}
It is easy to compute the second partial derivatives of $\hat{h}$
at $\mathbf{0}$ using \eqref{regpotential2} and the boundary conditions
$\left.\delta\right\vert _{\Gamma}=\,s$ and $\left.\hat{h}\right\vert _{\Gamma}=\,r$.
This gives, omitting the evaluation point $\mathbf{0}$, 
\begin{equation}
\begin{cases}
\hat{h}_{xx}=\,r^{\prime\prime}\left(0\right),\\
\hat{h}_{xy}=\frac{h_{xx}\,s\left(0\right)-(a-b\,\gamma)\,\,r^{\prime\prime}\left(0\right)}{(b-c\,\gamma)},\\
\hat{h}_{yy}=\frac{h_{xy}\,s\left(0\right)}{(b-c\,\gamma)}-\frac{(a-b\,\gamma)\,h_{xx}\,s\left(0\right)+(a-b\,\gamma)^{2}\,\,r^{\prime\prime}\left(0\right)}{(b-c\,\gamma)^{2}}. 
\end{cases}
\label{exex}
\end{equation}
Then, we must have
\begin{equation}
\hat{h}_{xx}=\,r^{\prime\prime}\left(0\right)>0\label{hes2}
\end{equation}
and 
\begin{align}
\hat{h}_{xx}\,\hat{h}_{yy}-\hat{h}_{xy}^{2}= & \frac{\,s\left(0\right)\,\,r^{\prime\prime}\left(0\right)}{(b-c\,\gamma)^{2}}\left(\vphantom{\frac{h_{xx}^{2}\,\,s\left(0\right)}{\,r^{\prime\prime}\left(0\right)}}(a-b\,\gamma)\,h_{xx}+(b-c\,\gamma)\,h_{xy}-\frac{h_{xx}^{2}\,\,s\left(0\right)}{\,r^{\prime\prime}\left(0\right)}\right)>0.\label{hes2'}
\end{align}
Accordingly, since \eqref{condpos} holds by hypothesis, in order
to ensure \eqref{hes} it is enough to take 
\begin{equation}
r^{\prime\prime}\left(0\right)>\frac{h_{xx}^{2}\left(\mathbf{0}\right)\,\,s\left(0\right)}{\left[(a-b\,\gamma)\,h_{xx}+(c-b\,\gamma)\,h_{xy}\right]\left(\mathbf{0}\right)}.\label{pr2}
\end{equation}
\end{proof}
\bigskip{}

The next lemma gives a necessary and sufficient condition, in terms
of $H$, for the existence of a function $\gamma$ fulfilling \eqref{condpos}.
The proof can be found in the Appendix.

\begin{lemma} \label{lem:gamma}There exists a function $\gamma$
such that \eqref{condpos} holds if and only if 
\begin{equation}
\left[b\,h_{xx}+c\,h_{xy}\right]\left(\mathbf{0}\right)\neq0\quad\text{or}\quad h_{xx}\left(\mathbf{0}\right)>0.\label{condpos2}
\end{equation}
Moreover, in such a case, $\gamma\left(\mathbf{0}\right)$ can be
chosen such that 
\begin{equation}
\left\{ \begin{array}{l}
\left\vert \gamma\left(\mathbf{0}\right)\right\vert >\left\vert \left[\frac{a\,h_{xx}+b\,h_{xy}}{b\,h_{xx}+c\,h_{xy}}\right]\left(\mathbf{0}\right)\right\vert \\
\left[\gamma\,\left(b\,h_{xx}+c\,h_{xy}\right)\right]\left(\mathbf{0}\right)<0
\end{array}\right.\label{c1}
\end{equation}
if $\left[b\,h_{xx}+c\,h_{xy}\right]\left(\mathbf{0}\right)\neq0$,
and using the following table
\begin{equation}
\begin{tabular}{|c|c|c|c|}
\hline   &  \ensuremath{h_{xy}\left(\mathbf{0}\right)=0}  &  \ensuremath{h_{xy}\left(\mathbf{0}\right)<0}  &  \ensuremath{h_{xy}\left(\mathbf{0}\right)>0}\\
\hline  \ensuremath{b\left(\mathbf{0}\right)=0}  &  any  &  \ensuremath{\gamma(\mathbf{0})>0}  &  \ensuremath{\gamma(\mathbf{0})<0}\\
\hline  \ensuremath{b\left(\mathbf{0}\right)>0}  &  \ensuremath{\gamma(\mathbf{0})<\frac{a\left(\mathbf{0}\right)}{b\left(\mathbf{0}\right)}}  &  \ensuremath{\frac{b\left(\mathbf{0}\right)}{c\left(\mathbf{0}\right)}<\gamma(\mathbf{0})<\frac{a\left(\mathbf{0}\right)}{b\left(\mathbf{0}\right)}}  &  \ensuremath{\gamma(\mathbf{0})<\text{min}\left(\frac{a\left(\mathbf{0}\right)}{b\left(\mathbf{0}\right)},\frac{b\left(\mathbf{0}\right)}{c\left(\mathbf{0}\right)}\right)}\\
\hline  \ensuremath{b\left(\mathbf{0}\right)<0}  &  \ensuremath{\gamma(\mathbf{0})>\frac{a\left(\mathbf{0}\right)}{b\left(\mathbf{0}\right)}}  &  \ensuremath{\gamma(\mathbf{0})>\text{max}\left(\frac{a\left(\mathbf{0}\right)}{b\left(\mathbf{0}\right)},\frac{b\left(\mathbf{0}\right)}{c\left(\mathbf{0}\right)}\right)}  &  \ensuremath{\frac{a\left(\mathbf{0}\right)}{b\left(\mathbf{0}\right)}<\gamma(\mathbf{0})<\frac{b\left(\mathbf{0}\right)}{c\left(\mathbf{0}\right)}} 
\\\hline \end{tabular}\label{c2}
\end{equation}
if $h_{xx}\left(\mathbf{0}\right)>0$. All of these conditions are compatible with \eqref{gbc}. 
\end{lemma}

Summarizing, if \eqref{condpos2} holds, in order to find a solution
$\left(\delta,\gamma,\hat{h}\right)$ of \eqref{regkinetic2}, \eqref{regpotential2}
and \eqref{pc}, it is enough to take $\gamma$ satisfying \eqref{gbc}
and also \eqref{c1} or \eqref{c2}, as explained in the last lemma.
Thus, we have proved the following.

\begin{theorem}\label{ida} Consider an underactuated simple Hamiltonian
system with two degrees of freedom and a set of adapted coordinates
related to it. Then if \eqref{condpos2} holds, the system is stabilizable
at $\mathbf{\bar{0}}$, i.e. there exists a state feedback controller
$u$, defined at least around $\mathbf{\bar{0}}$, such that the related
closed-loop system is stable at $\mathbf{\bar{0}}$. Moreover, such
a stability can be ensured by the existence of a simple Lyapunov function.
\end{theorem}

\subsection{A necessary condition for the existence of a simple Lyapunov function}

\label{necs} Using the same notation as above, suppose that an underactuated
simple Hamiltonian system with two degrees of freedom $\left(H,Y\right)$
can be stabilized at $\overline{\mathbf{0}}$, and that such stabilization
is ensured by the existence of a simple Lyapunov function. More precisely,
suppose that there exists a controller $u=\lambda Y$ and a simple
Lyapunov function $\hat{H}$, both of them defined at least around
$\mathbf{\bar{0}}$, ensuring the stability of the related closed-loop
system. Then, it was shown in \cite{lyap-shap} that $u$ and $\hat{H}$
must be given by the energy shaping method (or more precisely, by
the so-called \emph{simple CH method}). In particular, $\hat{H}$
must be locally given, in adapted coordinates, by a solution $\left(\delta,\gamma,\hat{h}\right)$
of \eqref{regkinetic2}, \eqref{regpotential2} and \eqref{pc} {[}and
by \eqref{HV1-1} and \eqref{invv} and some function $l>0${]}. We want
to show from this fact that \eqref{condpos2} must be satisfied. To
do that, let us consider two cases. %[\IEEEsetlabelwidth{g does not satisfy (30)}]

\begin{enumerate}
\item {[}$\gamma$ does not satisfy \eqref{gbc}{]} If $\gamma(\mathbf{0})=b(\mathbf{0})/c(\mathbf{0})$,
then {[}recall \eqref{pd}{]}
\[
\left[a-b\,\gamma\right](\mathbf{0})=a(\mathbf{0})-b(\mathbf{0})\,\frac{b(\mathbf{0})}{c(\mathbf{0})}=\frac{\Delta(\mathbf{0})}{c(\mathbf{0})}>0.
\]
On the other hand, if we differentiate \eqref{regpotential2} and
evaluate the result at $\mathbf{0}$, we obtain 
\[
\left(a(\mathbf{0})-b(\mathbf{0})\,\gamma(\mathbf{0})\right)\,\hat{h}_{xx}(\mathbf{0})=h_{xx}(\mathbf{0})\,\delta(\mathbf{0}).
\]
(Recall that $\mathbf{0}$ is critical for $\hat{h}$ and $h$). As
a consequence, using that $\delta(\mathbf{0})>0$ and $\hat{h}_{xx}(\mathbf{0})>0$,
\[
h_{xx}(\mathbf{0})=\frac{\left[a-b\,\gamma\right](\mathbf{0})\,\hat{h}_{xx}(\mathbf{0})}{\delta(\mathbf{0})}>0.
\]
In other words, condition \eqref{condpos2} must hold.
\item {[}$\gamma$ satisfies \eqref{gbc}{]} Let us call $V$ the neighborhood
of $\mathbf{0}$ where the functions $\delta,\gamma,\hat{h}$ are
defined. Define $r\left(x\right):=\hat{h}\left(x,0\right)$ and $s(x):=\delta(x,0)$
for all $x$ such that $\left(x,0\right)\in V$. It is clear that
the domain of the last functions is an open neighborhood of $0$.
Then, as we saw in the previous section, differentiating \eqref{regpotential2}
and evaluating at $\mathbf{0}$ (and using that $\mathbf{0}$ is critical
for $h$), we arrive at \eqref{exex}. Thus, the positivity conditions
\eqref{hes} for $\hat{h}$ can be studied in terms of \eqref{hes2}
and \eqref{hes2'}. From the latter, and from the fact that $s\left(0\right)>0$
(since $\delta$ must be positive), it easily follows that \eqref{condpos}
is a necessary condition. But according to Lemma \ref{lem:gamma},
this says again that condition \eqref{condpos2} must be satisfied. 
\end{enumerate}
Combining above discussion with Theorem \ref{ida}, we have the following
characterization.

\begin{theorem} \label{vuelta} Under the conditions of Theorem \ref{ida},
$\left(H,Y\right)$ is stabilizable at $\mathbf{\bar{0}}$, and such
stability can be ensured by the existence of a simple Lyapunov
function, if and only if \eqref{condpos2} holds.\end{theorem}

\section{Asymptotic stabilizability}

\label{asymptotic} In Reference \cite{chang2}, it was shown that
the condition\footnote{Actually, a weaker condition is considered there (see Theorem III.3).}
 $\left[b\,h_{xx}+c\,h_{xy}\right]\left(\mathbf{0}\right)\neq0$
also implies asymptotic stability (as previously affirmed in \cite{ham},
without a proof). In any other case, in the same reference, an additional
condition is proposed to ensure this kind of stability. We show in
the next subsection that no condition other than \eqref{condpos2}
is needed to this end.

\subsection{No additional assumptions are needed}

Let $\left(\delta,\gamma,\hat{h}\right)$ be a solution of \eqref{regkinetic2},
\eqref{regpotential2} and \eqref{pc} defined around $\mathbf{0}$,
with $\gamma$ satisfying \eqref{gbc} and \eqref{condpos} and with \emph{boundary conditions}
given by functions $s$ and $r$, as described in the proof of Lemma \ref{lem:gamma-1}. That is to say, $\delta$ and $\hat{h}$ must
satisfy 
\begin{equation}
\delta\left(x,0\right)=\,s\left(x\right)\ \ \ \text{and\ \ \ }\hat{h}\left(x,0\right)=r\left(x\right),\label{bcond}
\end{equation}
with $\,s$ and $r$ fulfilling \eqref{pr1} and \eqref{pr2}. To ensure the existence of such a solution, we only must ask that \eqref{condpos2} hold. Let
$\hat{H}$ be given by \eqref{HV1-1} and \eqref{invv}, i.e. 
\begin{align}
\hat{H}\left(x,y,p_{x},p_{y}\right)= & \left[p_{x}^{2}\,\left(\frac{\delta\left(x,y\right)}{l\left(x,y\right)}-\gamma^{2}\left(x,y\right)\right)+2\,\gamma\left(x,y\right)\,p_{x}p_{y}+p_{y}^{2}\right]\,l\left(x,y\right)\label{Vl}\\
&+\hat{h}\left(x,y\right)\nonumber,
\end{align}
for some positive function $l$. To write down an explicit expression
for the controller, we must choose a non-negative function $\mu$
fulfilling \eqref{muci}. To that effect, it suffices to take 
\begin{equation*}
\mu\left(x,y,p_{x},p_{y}\right)=\varkappa\,\left(\gamma\left(x,y\right)\,p_{x}+p_{y}\right)^{2}\,l^{2}\left(x,y\right),%\label{mu}
\end{equation*}
for some positive constant $\varkappa$.

\begin{remark}\label{LS} 
According to the results of Ref. \cite{lyap-shap},
the subset $\mu^{-1}\left(0\right)$, which in this case is given by
\begin{equation}
\mu^{-1}\left(0\right)=\left\{ \left(x,y,p_{x},p_{y}\right):p_{x}\,\gamma\left(x,y\right)+p_{y}=0\right\}, \label{mudos}
\end{equation}
is the \emph{LaSalle surface} related to the Lyapunov function $\hat{H}$
(see \cite{khalil}). Note also that $\mathbf{\bar{0}}\in\mu^{-1}\left(0\right)$. 
\end{remark}

%\begin{remark}
%It is not difficult to see that any choice for $\mu$ different from \eqref{mu} making \eqref{muci} smooth would contain the set \eqref{mudos}. In fact, as it will be clear later, the only values that are relevant for proving the asymptotic stability of the closed-loop system related to the controller defined by \eqref{Y,sol,mudos}; are those of $\partial \mu/\partial p_y$ on points of the form $p_y=-\gamma(x,y)p_x$.%, so with this choice we are selecting essentialy the corresponding term in a Taylor expansion in powers of $p_y$. 
%\end{remark}

With all these elements, the state feedback controller $u$ adopts
the form 
\begin{equation}
u=\left(0,0,0,\lambda\right),\label{ila}
\end{equation}
where $\lambda$ is locally given as {[}see \eqref{sol}{]} 
\begin{align}
\lambda\left(x,y,p_{x},p_{y}\right)= & -\varkappa\,\left(\gamma\left(x,y\right)\,p_{x}+p_{y}\right)\,l\left(x,y\right) -\frac{\left\{ \hat{H},H\right\} \left(x,y,p_{x},p_{y}\right)}{\left(\gamma\left(x,y\right)\,p_{x}+p_{y}\right)\,l\left(x,y\right)}.\label{lamu}
\end{align}
Our next step will be to prove that the
functions $\delta$, $\gamma$, $l$ and $\hat{h}$ can be chosen
such that the closed-loop system defined by $u$ is asymptotically
stable around $\overline{\mathbf{0}}$, i.e. the origin $\overline{\mathbf{0}}$
is an asymptotically stable equilibrium point of the vector field
$X=X_{H}+u$. More precisely, we are going to show, without any additional
assumption other than \eqref{condpos2}, that boundary conditions
$s$ and $r$ {[}see \eqref{bcond}{]} and an open subset $T$ containing
$\mathbf{\bar{0}}$ can be chosen in such a way that the largest $X$-invariant\footnote{Recall that, given a manifold $P$ and a vector field $X$ on $P$,
a subset $S\subset P$ is $X$\emph{-invariant} if every integral
curve of $X$ with initial condition in $S$ is contained in $S$.} submanifold of $\mathcal{S}_{0}:=\mu^{-1}\left(0\right)\cap T$ is
the singleton $\left\{ \mathbf{\bar{0}}\right\} $. Taking into account Remark \ref{LS}, this would imply,
\emph{via} the LaSalle's invariance principle, that $\mathbf{\bar{0}}$
is (locally) asymptotically stabilizable for $X$ (see \cite{khalil}).
The proof will be based on the next two lemmas (the proof of the first
one is easy to derive, so we omit it for brevity).

\begin{lemma} Given a manifold $P$, a vector field $X$ on $P$,
a critical point $\alpha_{0}$ of $X$, and a submanifold $\mathcal{S}_{0}\subset P$
containing $\alpha_{0}$, let us define\footnote{Given a manifold $P$, by $TP$ we are denoting, as usual, its tangent
bundle.}
\begin{equation}
\mathcal{S}_{n}:=\{\alpha\in\mathcal{S}_{n-1}:X(\alpha)\in T\mathcal{S}_{n-1}\},\ \ \ n\in\mathbb{N},\label{Sn}
\end{equation}
where we are assuming that each $\mathcal{S}_{n}$ is a submanifold
of $\mathcal{S}_{n-1}$. Then, the largest $X$-invariant subset $I$
of $\mathcal{S}_{0}$ satisfies 
\[
\displaystyle\left\{ \alpha_{0}\right\} \subset I\subset
\bigcap_{n\in\mathbb{N}}
\mathcal{S}_{n}.
\]
In particular, if $\mathcal{S}_{k}=\left\{ \alpha_{0}\right\} $ for
some $k\in\mathbb{N}$, then $I=\left\{ \alpha_{0}\right\} $. \end{lemma}

\begin{lemma}\label{impo} There exist boundary conditions $s$ and
$r$, a function $\gamma$ and an open subset $T\ni\mathbf{\bar{0}}$
such that {[}see \eqref{Sn}{]}: 
\begin{itemize}
\item the subset $\mathcal{S}_{1}$ corresponding to $\mathcal{S}_{0}=\mu^{-1}\left(0\right)\cap T$
is a submanifold of $\mathcal{S}_{0}$; 
\item $\mathcal{S}_{2}$ is a submanifold of $\mathcal{S}_{1}$; 
\item $\mathcal{S}_{3}=\left\{ \mathbf{\bar{0}}\right\} $. 
\end{itemize}
It is enough to take $s$ and $r$ such that {[}besides \eqref{pr1}
and \eqref{pr2}{]} 
\begin{equation}
\frac{\,s^{\prime}\left(0\right)}{\,s\left(0\right)}\neq-\left[\frac{2\,\left(b-c\,\gamma\right)}{\Delta}\left(b_{x}-\gamma\,c_{x}-\frac{B\,c}{2\,(b-c\,\gamma)}\right)\right]\left(\mathbf{0}\right),\label{l1}
\end{equation}
and choose $\gamma(\mathbf{0})$ according to \eqref{gbc}, \eqref{condpos}
and the additional restriction 
\begin{equation}
\gamma\left(\mathbf{0}\right)\neq\frac{\left(a\left(\mathbf{0}\right),b\left(\mathbf{0}\right)\right)\,\mathbb{M}\,\begin{pmatrix}b\left(\mathbf{0}\right)\\
c\left(\mathbf{0}\right)
\end{pmatrix}}{\left(b\left(\mathbf{0}\right),c\left(\mathbf{0}\right)\right)\,\mathbb{M}\,\begin{pmatrix}b\left(\mathbf{0}\right)\\
c\left(\mathbf{0}\right)
\end{pmatrix}},\label{l2}
\end{equation}
where $\mathbb{M}$ is a positive-definite matrix given by 
\begin{equation}
\mathbb{M}=\begin{bmatrix}\,r^{\prime\prime}\left(0\right) & \frac{h_{xx}\,s\left(0\right)-(a-b\,\gamma)\,r^{\prime\prime}\left(0\right)}{(b-c\,\gamma)}\\
\frac{h_{xx}\,s\left(0\right)-(a-b\,\gamma)\,r^{\prime\prime}\left(0\right)}{(b-c\,\gamma)} & \frac{h_{xy}\,s\left(0\right)\,(b-c\,\gamma)-(a-b\,\gamma)\,h_{xx}\,s\left(0\right)+(a-b\,\gamma)^{2}\,\,r^{\prime\prime}\left(0\right)}{(b-c\,\gamma)^{2}}
\end{bmatrix}.\label{M}
\end{equation}
\end{lemma} 
\begin{proof}
According to \eqref{mudos}, $\mu^{-1}\left(0\right)$ can be described
as the zero set of the function 
\[
\mathfrak{F}(x,y,p_{x},p_{y}):=\gamma(x,y)\,p_{x}+p_{y}.
\]
We shall proceed in three steps. 
\begin{enumerate}
\item Let us consider the subset $Z_{1}\subset\mu^{-1}\left(0\right)$ such
that $\mathfrak{F}_{\ast}(X)(x,y,p_{x},p_{y})=0$, where $\mathfrak{F}_{\ast}$ is the tangent map of $\mathfrak{F}$ and {[}see \eqref{XH}
and \eqref{ila}{]} 
\begin{align}
X=X_{H}+u= & \left(\frac{\partial H}{\partial p_{x}},\frac{\partial H}{\partial p_{y}},-\frac{\partial H}{\partial x},-\left(\frac{\partial H}{\partial y}-\lambda\right)\right).\label{xxhy}
\end{align}
That is to say, $Z_{1}$ is given by the equation 
\[
\frac{\partial H}{\partial p_{x}}\,\frac{\partial\mathfrak{F}}{\partial x}+\frac{\partial H}{\partial p_{y}}\,\frac{\partial\mathfrak{F}}{\partial y}-\frac{\partial H}{\partial x}\,\frac{\partial\mathfrak{F}}{\partial p_{x}}-\left(\frac{\partial H}{\partial y}-\lambda\right)\,\frac{\partial\mathfrak{F}}{\partial p_{y}}=0,
\]
or equivalently by 
\begin{equation}
\frac{\partial H}{\partial p_{x}}\,\gamma_{x}\,p_{x}+\frac{\partial H}{\partial p_{y}}\,\gamma_{y}\,p_{x}-\frac{\partial H}{\partial x}\,\gamma-\frac{\partial H}{\partial y}+\lambda=0,\label{zs}
\end{equation}
and $p_{y}=-\gamma(x,y)\,p_{x}$. Using the explicit forms of $H$
and $\hat{H}$, it is easy to see that, on $\mu^{-1}\left(0\right)$,
\begin{equation}
\left\{ \begin{array}{l}
\frac{\partial H}{\partial x}=\frac{1}{2}B\,p_{x}^{2}+h_{x},\ \ \frac{\partial H}{\partial y}=\frac{1}{2}C\,p_{x}^{2}+h_{y},\\
\frac{\partial H}{\partial p_{x}}=\left(a-b\,\gamma\right)\,p_{x},\ \ \frac{\partial H}{\partial p_{y}}=\left(b-c\,\gamma\right)\,p_{x},
\end{array}\right.\label{derh}
\end{equation}
where $B$ is given by \eqref{AB} and 
\[
C:= a_{y}-2\gamma b_{y}+\gamma^{2}c_{y}.
\]
On the other hand, and according to \eqref{lamu}, the values of $\lambda$
on points of the form $(x,y,p_{x},-\gamma\left(x,y\right)\,p_{x})\in\mu^{-1}\left(0\right)$,
are given by 
\begin{align*}
\lambda & =-\lim_{p_{y}\rightarrow-\gamma\,p_{x}}\frac{\left\{ \hat{H},H\right\} \left(x,y,p_{x},p_{y}\right)}{\left(\gamma\left(x,y\right)\,p_{x}+p_{y}\right)\,l\left(x,y\right)}\\
\\
 & =-\frac{1}{l\left(x,y\right)}\frac{\partial\{\hat{H},H\}}{\partial p_{y}}(x,y,p_{x},-\gamma\left(x,y\right)\,p_{x}).%\label{llim}
\end{align*}
So, by lengthy, but straightforward calculations, from \eqref{Vl},
\eqref{lamu} and \eqref{derh} we have that 
\begin{align}
\lambda= & \left[\frac{1}{2}(B\,\gamma+C)-\gamma_{x}\,\left(a-b\,\gamma\right)-\gamma_{y}\,\left(b-c\gamma\right)\right. \label{lat}\\
 & \left.+\frac{2\,(b_{x}-\gamma\,c_{x})\,\delta-b\,\delta_{x}-c\,\delta_{y}}{2\,l}\right]\,p_{x}^{2}+\gamma\,h_{x}+h_{y}-\frac{b\,\hat{h}_{x}+c\,\hat{h}_{y}}{\,l}.\nonumber
\end{align}
We are omitting, for simplicity, the evaluation point for the involved
functions. Finally, combining \eqref{zs}, \eqref{derh} and \eqref{lat},
we have at $\mu^{-1}\left(0\right)$ that
\begin{align*}
\mathfrak{F}_{\ast}(X)= & \frac{1}{l}\left[(b_{x}-\gamma\,c_{x})\,\delta-\frac{b\,\delta_{x}+c\,\delta_{y}}{2}\right]\,p_{x}^{2}-\frac{b\,\hat{h}_{x}+c\,\hat{h}_{y}}{\,l}.
\end{align*}
Thus, $\,Z_{1}$ is given by the equations 
\begin{equation}
\left\{ \begin{array}{l}
\gamma(x,y)\,p_{x}+p_{y}=0,\\
K(x,y)\,p_{x}^{2}-L(x,y)=0,
\end{array}\right.\label{kpl}
\end{equation}
where 
\begin{equation}
K=(b_{x}-\gamma\,c_{x})\,\delta-\frac{b\,\delta_{x}+c\,\delta_{y}}{2}\label{kl}
\end{equation}
and 
\begin{equation}
L=b\,\hat{h}_{x}+c\,\hat{h}_{y}.\label{kl'}
\end{equation}
In consequence, $\,Z_{1}$ can be defined by the zero set of the function
\[
\mathfrak{G}(x,y,p_{x},p_{y}):=\left(\mathfrak{G}_{1}(x,y,p_{x},p_{y}),\mathfrak{G}_{2}(x,y,p_{x},p_{y})\right).
\]
with 
\begin{align*}
\mathfrak{G}_{1}(x,y,p_{x},p_{y}) & :=\gamma(x,y)\,p_{x}+p_{y},\\
\mathfrak{G}_{2}(x,y,p_{x},p_{y}) & :=K(x,y)\,p_{x}^{2}-L(x,y).
\end{align*}
We want to see that its related tangent map (omitting the evaluation
point of the involved functions) 
\begin{equation}
\mathfrak{G}_{\ast}=\begin{pmatrix}\gamma_{x}\,p_{x} & \gamma_{y}\,p_{x} & \gamma & 1\\
K_{x}\,p_{x}^{2}-L_{x} & K_{y}\,p_{x}^{2}-L_{y} & 2\,K\,p_{x} & 0
\end{pmatrix}\label{gstar}
\end{equation}
has maximal rank around $\mathbf{\bar{0}}$. In that case, the Implicit
Function Theorem would ensure that this zero locus is the graph of
a smooth function (and hence a submanifold) when restricted to some
open neighborhood of $\mathbf{\bar{0}}$. Such tangent map is given
at $\mathbf{\bar{0}}$ by 
\begin{equation}
\mathfrak{G}_{\ast,\mathbf{\bar{0}}}=\begin{pmatrix}0 & 0 & \gamma(\mathbf{0}) & 1\\
-L_{x}(\mathbf{0}) & -L_{y}(\mathbf{0}) & 0 & 0
\end{pmatrix}.\label{g0}
\end{equation}
Note that the gradient of $L$ can be written 
\begin{equation}
\begin{pmatrix}L_{x}\\
L_{y}
\end{pmatrix}=\begin{pmatrix}\hat{h}_{xx} & \hat{h}_{xy}\\
\hat{h}_{xy} & \hat{h}_{yy}
\end{pmatrix}\,\begin{pmatrix}b\\
c
\end{pmatrix}+\begin{pmatrix}b_{x} & c_{x}\\
b_{y} & c_{y}
\end{pmatrix}\,\begin{pmatrix}\hat{h}_{x}\\
\hat{h}_{y}
\end{pmatrix}.\label{gL}
\end{equation}
Since $\mathbf{0}$ is critical for $\hat{h}$, then, at $\mathbf{0}$,
\[
\begin{pmatrix}L_{x}\\
L_{y}
\end{pmatrix}=\begin{pmatrix}\hat{h}_{xx} & \hat{h}_{xy}\\
\hat{h}_{xy} & \hat{h}_{yy}
\end{pmatrix}\,\begin{pmatrix}b\\
c
\end{pmatrix}.
\]
But we know that the Hessian matrix of $\hat{h}$ is positive-definite
and the function $c$ is always positive. So, the gradient of $L$
cannot vanish at $\mathbf{0}$. This implies that \eqref{g0} has
maximal rank at $\mathbf{\bar{0}}$. Consequently, there exists an
open subset $T_{1}$ containing $\mathbf{\bar{0}}$ such that $\,Z_{1}\cap T_{1}$
is a submanifold of $\mu^{-1}\left(0\right)\cap T_{1}$. 
\item Consider now the subset $\,Z_{2}\subset\,Z_{1}\cap T_{1}$ given by
$\mathfrak{G}_{\ast}(X)(x,y,p_{x},p_{y})=0$. Easy calculations show
that $\,Z_{2}$ is given by the points of $\,Z_{1}\cap T_{1}$ such
that 
\begin{align}
 & \left[\vphantom{\frac{1}{2}}(a-b\,\gamma)\,(K_{x}\,p_{x}^{2}-L_{x})+(b-c\,\gamma)\,(K_{y}\,p_{x}^{2}-L_{y})\right.\label{le}\\
 & \left.-2\,K\left(\frac{1}{2}B\,p_{x}^{2}+h_{x}\right)\right]\,p_{x}=0.\nonumber
\end{align}
We only need to evaluate the second row of $\mathfrak{G}_{\ast}$
{[}see \eqref{gstar}{]} on the components of $X$ {[}see \eqref{xxhy},
\eqref{derh} and \eqref{lat}{]}. In the following, we assume that
$K\left(\mathbf{0}\right)\neq0$. Observe that, since $\delta(x,0)=s(x)$,
we have $\delta_{x}(\mathbf{0})=s^{\prime}(0)$, and using \eqref{regpotential2}
at the origin 
\[
\delta_{y}\left(\mathbf{0}\right)=\left[\frac{B\,s(0)-(a-b\,\gamma)\,s^{\prime}\left(0\right)}{b-c\,\gamma}\right]\left(\mathbf{0}\right).
\]
So {[}see \eqref{kl}{]} 
\begin{align*}
K\left(\mathbf{0}\right)= & \left[\left(b_{x}-\gamma\,c_{x}-\frac{B\,c}{2\,(b-c\gamma)}\right)\,s\left(0\right)+\frac{1}{2}\left(\frac{\Delta}{b-c\gamma}\right)\,s^{\prime}\left(0\right)\right]\left(\mathbf{0}\right),
\end{align*}
and consequently, the condition $K\left(\mathbf{0}\right)\neq0$ is
equivalent to \eqref{l1}. Under such an assumption, we can replace
$p_{x}^{2}$ by $\frac{L}{K}$ in \eqref{le} {[}see \eqref{kpl}{]},
and we get 
\begin{align}
 & \left[(a-b\,\gamma)\,\left(K_{x}\,\frac{L}{K}-L_{x}\right)+(b-c\gamma)\,\left(K_{y}\,\frac{L}{K}-L_{y}\right)\right. \label{bed}\\
 & \left.-2\,K\left(\frac{1}{2}B\,\frac{L}{K}+h_{x}\right)\right]\,p_{x}=0,\nonumber
\end{align}
on some open subset containing $\mathbf{0}$ (where $K$ is non vanishing).
Moreover, since $L(\mathbf{0})=0$ {[}see \eqref{kl'}{]} and $h_{x}\left(\mathbf{0}\right)=0$,
we have at $\mathbf{0}$ that 
\[
\frac{L_{x}\,K-L\,K_{x}}{K^{2}}=\frac{L_{x}}{K},\ \ \frac{L_{y}\,K-L\,K_{y}}{K^{2}}=\frac{L_{y}}{K},
\]
and then, the bracketed expression in \eqref{bed} takes the following
form at $\mathbf{0}$ 
\[
(a-b\,\gamma)\,L_{x}+(b-c\,\gamma)\,L_{y}.
\]
Using \eqref{gL}, this in turn may be written as 
\[
\left(a,b\right)\,\begin{bmatrix}\hat{h}_{xx} & \hat{h}_{xy}\\
\hat{h}_{xy} & \hat{h}_{yy}
\end{bmatrix}\,\begin{pmatrix}b\\
c
\end{pmatrix}-\gamma\,\left(b,c\right)\,\begin{bmatrix}\hat{h}_{xx} & \hat{h}_{xy}\\
\hat{h}_{xy} & \hat{h}_{yy}
\end{bmatrix}\,\begin{pmatrix}b\\
c
\end{pmatrix}.
\]
Then, if we assume that, at $\mathbf{0}$, 
\[
\gamma\neq\frac{\left(a,b\right)\begin{bmatrix}\hat{h}_{xx} & \hat{h}_{xy}\\
\hat{h}_{xy} & \hat{h}_{yy}
\end{bmatrix}\begin{pmatrix}b\\
c
\end{pmatrix}}{\left(b,c\right)\begin{bmatrix}\hat{h}_{xx} & \hat{h}_{xy}\\
\hat{h}_{xy} & \hat{h}_{yy}
\end{bmatrix}\begin{pmatrix}b\\
c
\end{pmatrix}}
\]
it follows that \eqref{bed} will hold only if $p_{x}=0$ around $(x,y)=\mathbf{0}$.
It is worth mentioning that this condition is compatible with \eqref{gbc}
and \eqref{condpos}. Note that, using \eqref{exex}, the condition above
is given precisely by \eqref{l2} and \eqref{M}. In conclusion, there
exists an open neighborhood $T_{2}^{\prime}$ (which contains the
point $\mathbf{\bar{0}}$) such that the subset $\,Z_{2}\cap T_{2}^{\prime}$
is given by 
\[
\gamma\,p_{x}+p_{y}=0,\ \ K\,p_{x}^{2}-L=0,\ \ p_{x}=0,
\]
or equivalently 
\begin{equation}
p_{x}=p_{y}=L=0.\label{ppL}
\end{equation}
This means that $\,Z_{2}\cap T_{2}^{\prime}$ can be described as
the zero set of the function 
\[
\mathfrak{H}(x,y,p_{x},p_{y}):=\left(p_{y},p_{x},L(x,y)\right).
\]
The tangent map of $\mathfrak{H}$ at $\mathbf{\bar{0}}$ is given
by 
\[
\mathfrak{H}_{\ast,\mathbf{\bar{0}}}=\begin{pmatrix}0 & 0 & 0 & 1\\
0 & 0 & 1 & 0\\
L_{x}(\mathbf{0}) & L_{y}(\mathbf{0}) & 0 & 0
\end{pmatrix}.
\]
Again, since $L_{x}(\mathbf{0})$ and $L_{y}(\mathbf{0})$ cannot
be both zero, we conclude that $\mathfrak{H}_{\ast,\mathbf{\bar{0}}}$
has maximal rank. Thus, there exists inside $T_{2}^{\prime}$ an open
neighborhood $T_{2}$ of $\mathbf{\bar{0}}$ such that $Z_{2}\cap T_{2}$
is a submanifold of $\,Z_{1}\cap T_{1}\cap T_{2}$. 
\item Now, consider the subset $\,Z_{3}\subset\,Z_{2}\cap T_{2}$ defined
by $\mathfrak{H}_{\ast}(X)(x,y,p_{x},p_{y})=0$. Using \eqref{xxhy},
\eqref{derh} and \eqref{lat}, it follows that, along $\,Z_{2}\cap T_{2}$
{[}see \eqref{ppL}{]} 
\[
X=\left(0,0,-h_{x},\gamma\,h_{x}\right),
\]
so, in order for $\mathfrak{H}_{\ast}(X)$ to vanish, it is necessary
that $h_{x}=0$. But, if this is the case, using the potential matching
condition 
\[
(a-b\,\gamma)\,\hat{h}_{x}+(b-c\,\gamma)\,\hat{h}_{y}=\delta\,h_{x},
\]
or equivalently, $a\,\hat{h}_{x}+b\,\hat{h}_{y}-\gamma\,L=\delta\,h_{x}$,
we have on $\,Z_{3}$ that 
\[
L=b\,\hat{h}_{x}+c\,\hat{h}_{y}=0\ \ \text{and\ \ }a\,\hat{h}_{x}+b\,\hat{h}_{y}=0,
\]
i.e. 
\[
\mathbb{H}\,\left(\begin{array}{l}
\hat{h}_{x}\\
\hat{h}_{y}
\end{array}\right)=0.
\]
Calling $\pi$ the projection of $\mathbb{R}^{4}$ onto the first
two components, we can say that above identity holds if and only if
all the points of $\pi(Z_{3})$ are critical for $\hat{h}$. By the
\emph{Morse Lemma}, since $\mathbf{0}$ is a non-degenerate critical
point of $\hat{h}$ [recall \eqref{pc}], there exists a neighborhood $V$ of $\mathbf{0}$
such that $\pi\left(\,Z_{3}\right)\cap V=\left\{ \mathbf{0}\right\} $.
But $p_{x}=p_{y}=0$ on $\,Z_{3}$, which implies that $Z_{3}\cap T_{3}=\left\{ \mathbf{\bar{0}}\right\} $
for $T_{3}:=\pi^{-1}\left(V\right)$. 
\end{enumerate}
Summing up, if we define $T:=T_{1}\cap T_{2}\cap T_{3}$ and $\mathcal{S}_{0}:=\mu^{-1}\left(0\right)\cap T$,
from \eqref{Sn} we obtain $\mathcal{S}_{1}=Z_{1}\cap T$, which is
a submanifold of $\mathcal{S}_{0}$, $\mathcal{S}_{2}=Z_{2}\cap T$,
which is a submanifold of $\mathcal{S}_{1}$, and $\mathcal{S}_{3}=\left\{ \mathbf{\bar{0}}\right\} $.
Hence, the three points of the lemma follow. 
\end{proof}

Concluding, if \eqref{condpos2} holds, asymptotic stabilizability
is ensured. Reciprocally, if we can ensure asymptotic stabilizability
by the existence of a simple Lyapunov function, then we can also ensure
stabilizability, and Theorem \ref{vuelta} implies that \eqref{condpos2}
holds. In other terms,

\begin{theorem} Under the conditions of Theorem \ref{ida}, $\left(H,Y\right)$
is asymptotically stabilizable at $\mathbf{\bar{0}}$, and such a
stability can be ensured by the existence of a simple Lyapunov function,
if and only if $H$ satisfies \eqref{condpos2}. \end{theorem}

\subsection{Example: the inertia wheel pendulum}

Now, we illustrate our results with a concrete underactuated system
$\left(H,Y\right)$ with two degrees of freedom, the \emph{inertia
wheel pendulum}: 
\begin{itemize}
\item the configuration space is $Q=S^{1}\times S^{1}$, whose natural almost-global
coordinates will be denoted $\left(\theta,\psi\right)$; 
\item the Hamiltonian is 
\begin{align*}
H\left(\theta,\psi,p_{\theta},p_{\psi}\right)= & \frac{1}{2}\left(p_{\theta},p_{\psi}\right)\,\begin{bmatrix}a & b\\
b & c
\end{bmatrix}\,\begin{pmatrix}p_{\theta}\\
p_{\psi}
\end{pmatrix} +M\,(1+\cos\theta),
\end{align*}
where $a,b,c,M$ are constants and $a,b,M,ac-b^{2}$ are strictly
greater than zero; 
\item and the space of actuators is given by the subbundle spanned by the
vector field $Y$ taking the constant value $\left(0,0,0,1\right)$. 
\end{itemize}
We shall find, by using the energy shaping method, a state feedback
controller $u$ for this system and a related simple Lyapunov function
$\hat{H}$ which make the closed-loop system $X_{H}+u$ asymptotically
stable at $\left(\theta,\psi,p_{\theta},p_{\psi}\right)=\left(0,0,0,0\right)=\mathbf{\bar{0}}$.

Replacing $x$ by $\theta$ and $y$ by $\psi$, \eqref{condpos2}
in this case says that (because $h_{\theta\psi}\left(\mathbf{0}\right)=0$)
\[
b\,h_{\theta\theta}\left(\mathbf{0}\right)\neq0\quad\text{or}\quad h_{\theta\theta}\left(\mathbf{0}\right)>0,
\]
which is equivalent to $h_{\theta\theta}\left(\mathbf{0}\right)\neq0$,
since $b\neq0$. And it does hold, because $h_{\theta\theta}\left(\mathbf{0}\right)=-M\neq0$.
Then, the inertia wheel pendulum can be asymptotically stabilized
around $\mathbf{\bar{0}}$, as it is well known. On the other hand,
according to \eqref{AB}, we have that $B=0$. So, the kinetic and
potential matching conditions read {[}see \eqref{regkinetic2} and
\eqref{regpotential2}{]} 
\begin{equation}
(a-b\,\gamma)\,\delta_{\theta}+(b-c\,\gamma)\,\delta_{\psi}=0,\label{kef}
\end{equation}
and 
\begin{equation}
(a-b\,\gamma)\,\hat{h}_{\theta}+(b-c\,\gamma)\,\hat{h}_{\psi}=-M\,\delta\,\sin\theta,\label{pef}
\end{equation}
respectively. Let us construct a solution $\left(\delta,\gamma,\hat{h}\right)$
of above equations, with $\delta>0$ and $\hat{h}$ positive-definite
w.r.t. $\mathbf{0}$. We shall take $\gamma$ constant. Following
the steps of \S\ref{exandun}, it is enough to take $\gamma$
such that $\gamma\neq b/c$ {[}see \eqref{gbc}{]} and, using \eqref{c1}
of Lemma \ref{lem:gamma} {[}since $h_{\theta\theta}\left(\mathbf{0}\right)=-M<0${]},
also ask that 
\[
\left\vert \gamma\right\vert >\frac{a}{b\,}\ \ \ \ \ \text{and\ \ \ }-\gamma\,b\,M<0.
\]
The second inequality says that $\gamma$ is positive, so, above equations
only impose the condition $\gamma>a/b$. Note also that, since $a,b,c,ac-b^{2}>0$,
we have that $a/b>b/c$. Hence, all the conditions on $\gamma$ reduce
to 
\[
\gamma>\frac{a}{b\,}.
\]

\begin{remark} Note that, for this system, we cannot take $\gamma=b/c$.
In fact, in such a case, according to the calculations we made in \S\ref{necs}, the positivity of $\delta$ and $\hat{h}$ would
impose that $h_{\theta\theta}\left(\mathbf{0}\right)>0$, which is
not true. \end{remark}

Regarding the boundary conditions defining $\delta$ and $\hat{h}$,
i.e. the functions $s$ and $r$, respectively, we must ask {[}see
\eqref{pr1} and \eqref{pr2}{]} 
\begin{equation}
\left\{ \begin{array}{l}
s\left(0\right)>0,\\
r\left(0\right)=\,r^{\prime}(0)=0\qquad\text{and}\\
r^{\prime\prime}\left(0\right)>\frac{h_{\theta\theta}^{2}\left(\mathbf{0}\right)\,s\left(0\right)}{(a-b\,\gamma)\,h_{\theta\theta}\left(\mathbf{0}\right)}=-\frac{M\,s\left(0\right)}{a-b\,\gamma}
\end{array}\right.\label{f1}
\end{equation}
And to ensure asymptotic stabilizability, according to Eqs. \eqref{l1}, \eqref{l2} and \eqref{M} of Lemma \ref{impo}, we ask that
\[
\frac{s^{\prime}\left(0\right)}{s\left(0\right)}\neq0,\ \ \ \text{i.e.\ \ }s^{\prime}\left(0\right)\neq0,
\]
and that $\gamma$ satisfy
\begin{equation}
\gamma\neq\frac{\eta^{2}\,ab\,r^{\prime\prime}\left(0\right)-\eta\,\left[M\,s\left(0\right)+\zeta\,r^{\prime\prime}\left(0\right)\right]\,\left(ac+b^{2}\right)+\zeta\,M\,s\left(0\right)+\zeta^{2}\,\,r^{\prime\prime}\left(0\right)\,bc}{\eta^{2}\,b^{2}\,r^{\prime\prime}\left(0\right)-2\,\eta\,\left[M\,s\left(0\right)+\zeta\,r^{\prime\prime}\left(0\right)\right]\,bc+\zeta\,M\,s\left(0\right)+\zeta^{2}\,\,r^{\prime\prime}\left(0\right)\,c^{2}},\label{f2}
\end{equation}
where $\zeta:=a-b\,\gamma$ and $\eta:=b-c\,\gamma$.
Thus, take any number $\gamma>a/b$, any function $s$ such that $\,s\left(0\right)>0$
and $s^{\prime}\left(0\right)\neq0$, and any function\footnote{Additionally, the functions $s$ and $r$ may be taken with period
$2\pi$ in order to look for a quasi-global solution.} $r$ such that $r\left(0\right)=\,r^{\prime}(0)=0$ and $r^{\prime\prime}\left(0\right)$
satisfying \eqref{f1} and \eqref{f2}, and let us apply the Method
of Characteristics to \eqref{kef} and \eqref{pef}, with boundary
conditions on $\psi=0$ given by $s$ and $r$. The characteristic
equations for \eqref{kef} are 
\[
\begin{cases}
\dot{\theta}=a-b\,\gamma, & \theta(0)=\theta_{0},\\
\dot{\psi}=b-c\,\gamma, & \psi(0)=0,\\
\dot{\delta}=0, & \delta(0)=s(\theta_{0}).
\end{cases}
\]
Then 
\[
\theta(t)=(a-b\,\gamma)\,t+\theta_{0},\ \ \ \psi(t)=(b-c\,\gamma)\,t,
\]
and defining $\Upsilon:=\left(a-b\,\gamma\right)/\left(b-c\,\gamma\right)$
we find 
\begin{equation}
\delta(\theta,\psi)=s(\theta-\Upsilon\,\psi).\label{delf}
\end{equation}
The characteristic equation for \eqref{pef} (and for $\delta$ given
above) is 
\[
\dot{\hat{h}}=-M\,s(\theta_{0})\,\sin\left((a-b\,\gamma)\,t+\theta_{0}\right),
\]
and integrating we obtain 
\begin{align}
\hat{h}(\theta,\psi)= & \frac{M\,s(\theta-\Upsilon\,\psi)}{a-b\,\gamma}\left(\cos\theta-\cos(\theta-\Upsilon\,\psi)\right)+\,r(\theta-\Upsilon\,\psi).\label{vf}
\end{align}

Finally, with $\delta$ and $\hat{h}$ given by \eqref{delf} and
\eqref{vf}, and considering any positive function $l$, we have from
\eqref{Vl}, \eqref{ila} and \eqref{lamu} the controller $u$ and
the Lyapunov function $\hat{H}$ we are looking for.

\appendix

\section{Proof of Lemma \ref{lem:gamma}}

Suppose first that $h_{xx}\left(\mathbf{0}\right)\leq0$ and $\left(b\,h_{xx}+c\,h_{xy}\right)\left(\mathbf{0}\right)=0$.
Then, omitting the evaluation point $\mathbf{0}$, 
\begin{align*}
(a-b\,\gamma)\,h_{xx}+(b-c\,\gamma)\,h_{xy} & =a\,h_{xx}+b\,h_{xy}-\gamma\,(b\,h_{xx}+c\,h_{xy})\\
 & =\left(a-\frac{b^{2}}{c}\right)\,h_{xx}=\frac{\Delta}{c}\,h_{xx}\leq0
\end{align*}
for any function $\gamma$. This proves the first implication of the
lemma (by denying the second one). For the converse, suppose first
that $b\,h_{xx}+c\,h_{xy}\neq0$. Since 
\[
(a-b\,\gamma)\,h_{xx}+(b-c\,\gamma)\,h_{xy}=a\,h_{xx}+b\,h_{xy}-\gamma\,(b\,h_{xx}+c\,h_{xy}),
\]
for any $\gamma$ with sign opposite to $b\,h_{xx}+c\,h_{xy}$ and
such that 
\[
|\gamma|>\left\vert \frac{a\,h_{xx}+b\,h_{xy}}{b\,h_{xx}+c\,h_{xy}}\right\vert ,
\]
we have that \eqref{condpos} holds. This implies \eqref{c1}. Now,
suppose that $h_{xx}>0$. If $h_{xy}=0$, since $a>0$, it is clear
that it is enough to choose $\gamma$ such that $a>b\gamma$. If instead
$h_{xy}\neq0$, we distinguish three cases: $b=0$, $b>0$ and $b<0$. 
\begin{itemize}
\item If $b=0$, then $a-b\gamma=a>0$ and $b-c\gamma=-c\gamma$, and consequently
it is sufficient to choose $\gamma$ with opposite sign to $h_{xy}$. 
\item If $b>0$, we show that it is possible to take $\gamma$ so as to
fulfill one of the following expressions 
\[
a-b\gamma>0\quad\text{ and }\quad b-c\gamma>0,
\]
or 
\[
a-b\gamma>0\quad\text{ and }\quad b-c\gamma<0.
\]
In order to make $a-b\gamma>0$, we need $\gamma<\frac{a}{b}$. If
in addition $b-c\gamma>0$, then $\gamma<\frac{b}{c}$. Hence, it
suffices to take $\gamma<\text{min}\left(\frac{a}{b},\frac{b}{c}\right)$.
On the contrary, if $b-c\gamma<0$, then we can take $\frac{b}{c}<\gamma<\frac{a}{b}$,
which is always possible because 
\[
ac-b^{2}>0\;\text{ and }\;b>0\qquad\Rightarrow\qquad\frac{b}{c}<\frac{a}{b}.
\]
Thus, if $h_{xy}>0$ we choose $\gamma<\text{min}\left(\frac{a}{b},\frac{b}{c}\right)$
and if $h_{xy}<0$ we take $\frac{b}{c}<\gamma<\frac{a}{b}$. In both
cases we get the desired result. 
\item If $b<0$, then $a-b\gamma>0$ implies $\gamma>\frac{a}{b}$. If in
addition $b-c\gamma>0$, then $\gamma<\frac{b}{c}$ and so we can
take $\gamma$ such that $\frac{a}{b}<\gamma<\frac{b}{c}$, which
is always possible because 
\[
ac-b^{2}>0\;\text{ and }\;b<0\qquad\Rightarrow\qquad\frac{a}{b}<\gamma<\frac{b}{c}.
\]
On the contrary, if $b-c\gamma<0$, then $\gamma>\frac{b}{c}$ and
it is sufficient to choose $\gamma>\text{max}\left(\frac{a}{b},\frac{b}{c}\right)$.
Again, both cases lead to the desired result. All these alternatives
give the table \eqref{c2}. 
\end{itemize}
The last assertion of the lemma is immediate, because all conditions
on $\gamma$ are inequalities.

%%%%% Bibliography %%%%%%%%%%%%%%%%%%%%%%%%%

\end{document}